\input amstex
\input amsppt.sty
\magnification=\magstep1
\hsize=30truecc
\vsize=22.2truecm
\baselineskip=16truept
\NoBlackBoxes
\TagsOnRight \pageno=1 \nologo
\def\Z{\Bbb Z}
\def\N{\Bbb N}

\def\C{\Bbb C}
\def\l{\left}
\def\r{\right}
\def\bg{\bigg}
\def\({\bg(}
\def\[{\bg\lfloor}
\def\){\bg)}
\def\]{\bg\rfloor}
\def\t{\text}
\def\f{\frac}

\def\bi{\binom}
\def\eq{\equiv}

\def\ls{\leqslant}
\def\gs{\geqslant}
\def\mo{\roman{mod}}

\def\al{\alpha}
\def\da{\delta}

\def\m#1#2{\thickfracwithdelims\{\}\thickness0{#1}{#2}_m}

\def\Proof{\noindent{\it Proof}}

\def\Remark{\medskip\noindent{\it  Remark}}

\hbox {Preprint, {\tt arXiv:0909.3808}}
\bigskip
\topmatter
\title Various congruences involving binomial coefficients and higher-order Catalan numbers\endtitle
\author Zhi-Wei Sun\endauthor
\leftheadtext{Zhi-Wei Sun}
\rightheadtext{Binomial coefficients and higher-order Catalan numbers}
\affil Department of Mathematics, Nanjing University\\
 Nanjing 210093, People's Republic of China
  \\  zwsun\@nju.edu.cn
  \\ {\tt http://math.nju.edu.cn/$\sim$zwsun}
\endaffil
\abstract Let $p$ be a prime and let $a$ be a positive integer.
In this paper we investigate $\sum_{k=0}^{p^a-1}\bi{(h+1)k}{k+d}/m^k$ modulo a prime $p$,
where $d$ and $m$ are integers with $-h<d\ls p^a$ and $m\not\eq0\ (\mo\ p)$.
We also study congruences involving higher-order Catalan numbers
$C_k^{(h)}=\f1{hk+1}\bi{(h+1)k}k$ and $\bar C_k^{(h)}=\frac h{k+1}\bi{(h+1)k}k$ .
Our tools include linear recurrences and the theory of cubic residues.
Here are some typical results in the paper.
(i) If $p^a\eq1\ (\mo\ 6)$ then
$$\sum_{k=0}^{p^a-1}\f{\bi{3k}k}{6^k}\eq 2^{(p^a-1)/3}\ (\mo\ p)
\ \ \t{and}\ \ \sum_{k=1}^{p^a-1}\f{\bar C_k^{(2)}}{6^k}\eq0\ (\mo\ p).$$
Also, $$\sum_{k=0}^{p^a-1}\f{\bi{3k}k}{7^k}\eq\cases-2&\t{if}\ p^a\eq\pm2\ (\mo\ 7),
\\1&\t{otherwise}.\endcases$$
(ii) We have
$$\sum_{k=0}^{p^a-1}\f{\bi{4k}k}{5^k}\eq\cases1\ (\mo\ p)&\t{if}\ p\not=11\ \t{and}\ p^a\eq1\ (\mo\ 5),
\\-1/11\ (\mo\ p)&\t{if}\ p^a\eq2,3\ (\mo\ 5),
\\-9/11\ (\mo\ p)&\t{if}\ p^a\eq4\ (\mo\ 5).
\endcases$$
Also, $$\sum_{k=0}^{p^a-1}\f{C_k^{(3)}}{5^k}\eq\cases1\ (\mo\ p)&\t{if}\ p^a\eq1,3\ (\mo\ 5),
\\-2\ (\mo\ p)&\t{if}\ p^a\eq2\ (\mo\ 5),
\\0\ (\mo\ p)&\t{if}\ p^a\eq4\ (\mo\ 5).\endcases
$$

\endabstract
\thanks 2010 {\it Mathematics Subject Classification}.\,Primary 11B65;
Secondary 05A10,\,11A07.
\newline\indent Supported by the National Natural Science
Foundation (grant 10871087) and the Overseas Cooperation Fund (grant 10928101) of China.
\endthanks
\endtopmatter
\document

\heading{1. Introduction}\endheading

Let $p$ be a prime. Via a sophisticated combinatorial identity, H. Pan and Z. W. Sun [PS]
proved that
$$\sum_{k=0}^{p-1}\bi{2k}{k+d}\eq\l(\f{p-d}3\r)\ (\mo\ p)\quad \t{for}\ d=0,\ldots,p,$$
where $(-)$ is the Jacobi symbol.
Let $a\in\Z^+=\{1,2,3,\ldots\}$ and $d\in\{0,\ldots,p^a\}$. Recently Sun and R. Tauraso [ST1] used a new approach
to determine $\sum_{k=0}^{p^a-1}\bi{2k}{k+d}$ mod $p^2$; they [ST2]
also studied $\sum_{k=1}^{p^a-1}\bi{2k}{k+d}/m^k$ modulo $p$
via Lucase sequences, where $m$ is an integer not divisible by $p$.

Quite recently, L. Zhao, Pan and Sun [ZPS] proved that if $p\not=2,5$ is a prime then
$$\sum_{k=1}^{p-1}2^k\bi{3k}k\eq\f65\l(\l(\f{-1}p\r)-1\r)\ (\mo\ p)$$
and
$$\sum_{k=1}^{p-1}2^{k-1}C_k^{(2)}\eq\l(\f{-1}p\r)-1\ (\mo\ p),$$
where $C_k^{(2)}=\bi{3k}k/(2k+1)\ (k\in\N=\{0,1,2,\ldots\})$ are Catalan numbers of order 2.

In general, (the first-kind) Catalan numbers of order $h\in\Z^+$ are given by
$$C_k^{(h)}=\f1{hk+1}\bi{(h+1)k}k=\bi{(h+1)k}k-h\bi{(h+1)k}{k-1}\ \ (k\in\N).$$
(As usual, $\bi x{-n}=0$ for $n=1,2,\ldots$.)
We also define the second-kind Catalan numbers of order $h$ as follows:
$$\bar C_k^{(h)}=\f h{k+1}\bi{(h+1)k}k=h\bi{(h+1)k}k-\bi{(h+1)k}{k+1}\ \ (k\in\N).$$
Those $C_k=C_k^{(1)}=\bar C_k^{(1)}$ are ordinary Catalan numbers which have lots of combinatorial interpretations
(see, e.g., Stanley [St]).

Let $p$ be a prime and $a$ a positive integer.
In this paper we mainly investigate $\sum_{k=0}^{p^a-1}\bi{3k}k/m^k$ mod $p$
for all $m\in\Z$ with $m\not\eq0\ (\mo\ p)$,
and determine $\sum_{k=0}^{p^a-1}\bi{4k}k/5^k$ and $\sum_{k=0}^{p^a-1}C_k^{(3)}/5^k$ modulo $p$.
Our approach involves third-order and fourth order recurrences and the theory of cubic residues.

Now we introduce some basic notations throughout this paper. For a positive integer $n$, we use $\Z_n$ to denote the set
of all rational numbers whose denominators are relatively prime to $n$. Thus, if $p$ is a prime then $\Z_p$ is the ring of rational
$p$-adic integers. For a predicate $P$, we let
$$[P]=\cases1&\t{if}\ P\ \t{holds},\\0&\t{otherwise}.\endcases$$
Thus $[m=n]$ coincides with the Kronecker $\da_{m,n}$.

Our first theorem is a further extension of the above-mentioned congruences of Zhao, Pan and Sun.

\proclaim{Theorem 1.1} Let $p$ be an odd prime and let $a\in\Z^+$.  Let $c\in\Z_p$ with $c\not\eq0,-1,2\ (\mo\ p)$,
and set $c'=3/(2(c+1)(c-2))$. Then
$$\sum_{k=1}^{p^a-1}\f{c^{2k}}{(c+1)^{3k}}\bi{3k}k\eq c'\(1-\l(\f{4c+1}{p^a}\r)\) \ (\mo\ p),$$
$$\sum_{k=1}^{p^a-1}\f{c^{2k+1}}{(c+1)^{3k}}\bi{3k}{k-1}\eq(c'+1)\(1-\l(\f{4c+1}{p^a}\r)\)\ (\mo\ p),$$
$$\sum_{k=1}^{p^a-1}\f{c^{2k+2}}{(c+1)^{3k}}\bi{3k}{k+1}
\eq(c'(3c+2)+1)\(1-\l(\f{4c+1}{p^a}\r)\) \ (\mo\ p),
$$
and
$$\sum_{k=0}^{p^a-1}\f{c^{2k}}{(c+1)^{3k}}\bi{3k}{k+p^a}\eq
cc'\(\l(\f{4c+1}{p^a}\r)-1\)\ (\mo\ p).$$
\endproclaim
\Remark\ 1.1. Note that if $c=-1/4$ then $c^2/(c+1)^3=2^2/(2+1)^3$.

\medskip
Clearly Theorem 1.1 in the case $c=-1/2$ yields the two congruences  of Zhao, Pan and Sun [ZPS] mentioned above.
Applying Theorem 1.1 with $c=1,-2$ we obtain the following consequence.

\proclaim{Corollary 1.1} Let $p$ be an odd prime and let $a\in\Z^+$. Then
$$\align\sum_{k=1}^{p^a-1}\f{\bi{3k}k}{8^k}\eq&\f34\(\l(\f{p^a}5\r)-1\)\ (\mo\ p),
\\\sum_{k=1}^{p^a-1}\f{C_k^{(2)}}{8^k}\eq&\f54\(\l(\f{p^a}5\r)-1\)\ (\mo\ p),
\\\sum_{k=1}^{p^a-1}(-4)^k\bi{3k}k\eq&\f38\(1-\l(\f{p^a}7\r)\)\ (\mo\ p),
\\\sum_{k=1}^{p^a-1}(-4)^kC_k^{(2)}\eq&\f74\(1-\l(\f{p^a}7\r)\)\ (\mo\ p).
\endalign$$
\endproclaim

For a polynomial
$$f(x)=x^n+a_1x^{n-1}+\cdots+a_n=\prod_{i=1}^n(x-\al_i)\in\C[x],$$
its discriminant is defined by
$$D(f)=\prod_{1\ls i<j\ls n}(\al_i-\al_j)^2.$$
By Vitae's theorem and the fundamental theorem of symmetric polynomials, we can express $D(f)$ as a rational expression
involving the coefficients $a_1,\ldots,a_n$. For example, it is known that
$$D(x^3+a_1x^2+a_2x+a_3)=a_1^2a_2^2-4a_2^3-4a_1^3a_3-27a_3^2+18a_1a_2a_3.$$
If $f(x)=x^n+a_1x^{n-1}+\cdots+a_n\in\Z[x]$ and $p$ is an odd prime not dividing $D(f)$, then
$$\l(\f{D(f)}p\r)=(-1)^{n-r}$$
by Stickelberger's theorem (cf. [C]), where $r$ is the total number of monic irreducible factors of $f(x)$ modulo $p$.

Let $p$ be an odd prime and $m$ an integer with $m\not\eq0,27/4\ (\mo\ p)$.
Then $D=D((x+1)^3-3x^2)=(4m-27)m^2\not\eq0\ (\mo\ p)$. Suppose that there is no $c\in\Z_p$ such that
$mc^2\eq(c+1)^3\ (\mo\ p)$. Then the polynomial $(1+x)^3-mx^2$ is irreducible modulo $p$,
hence by the Stickelberger theorem we have $(\f Dp)=(-1)^{3-1}=1$. Thus $(\f{4m-27}p)=1$,
and hence $4m-27\eq (2t+1)^2\ (\mo\ p)$ for some $t\in\Z$. Note that $m\eq t^2+t+7\ (\mo\ p)$.

The following theorem deals with the case $m=6$ and $(\f{4m-27}{p})=1$.

\proclaim{Theorem 1.2} Let $p>3$ be a prime and let $a\in\Z^+$. Suppose that $p^a\eq 1\ (\mo\ 6)$. Then
$$\sum_{k=1}^{p^a-1}\f{\bi{3k}k}{6^k(k+1)}\eq\sum_{k=1}^{p^a-1}\f{\bi{3k}{k-1}}{6^k}\eq0\ (\mo\ p)$$
and
$$\sum_{k=1}^{p^a-1}\f{\bi{3k}k}{6^k}\eq2^{(p^a-1)/3}-1\eq\f12\sum_{k=1}^{p^a-1}\f{\bi{3k}{k+1}}{6^k}\ (\mo\ p).$$
\endproclaim

Now we need to introduce another notation.
For a positive integer $n\not\eq0\ (\mo\ 3)$ and $i\in\{0,1,2\}$, Z.-H. Sun [S98] investigated
$$C_i(n)=\bg\{k\in\Z_n:\ \l(\f{k+1+2\omega}n\r)_3=\omega^i\bg\},$$
where $\omega$ is the primitive cubic root $(-1+\sqrt{-3})/2$ of unity, and
$(\f{\cdot}n)_3$ is the cubic Jacobi symbol.
(The reader is referred to Chapter 9 of [IR, pp.\,108-137] for the basic theory of cubic residues.)
By [S98],  $k\in C_2(n)$ if and only if $-k\in C_1(n)$; also
$$C_0(n)\cup C_1(n)\cup C_2(n)=\{k\in\Z_n:\ k^2+3\ \t{is relatively prime to} \ n\}.$$

\proclaim{Theorem 1.3} Let $p>3$ be a prime and let $a\in\Z^+$. Let $m,t\in\Z_p$ with $t\not\eq -1/2\ (\mo\ p)$
and $m\eq t^2+t+7\not\eq0,6\ (\mo\ p)$.
Then
$$c=\f{2m^2-18m+27}{6t+3}\in C_0(p^a)\cup C_1(p^a)\cup C_2(p^a).$$
If $c\in C_0(p^a)$, then
$$\sum_{k=1}^{p^a-1}\f{\bi{3k}{k+d}}{m^k}\eq0\ (\mo\ p)\quad\t{for}\ d\in\{0,\pm1\},$$
and hence
$$\sum_{k=1}^{p^a-1}\f{C_k^{(2)}}{m^k}\eq\sum_{k=1}^{p^a-1}\f{\bar C_k^{(2)}}{m^k}\eq 0\ (\mo\ p).$$
When $\pm c\in C_1(p^a)$, we have
$$\sum_{k=1}^{p^a-1}\f{\bi{3k}{k+d}}{m^k}\eq\cases(\pm 3/(2t+1)-3)/2\ (\mo\ p)&\t{if}\ d=0,
\\\pm(m-6)/(2t+1)\ (\mo\ p)&\t{if}\ d=-1,
\\\pm 3/(2t+1)+3-m\ (\mo\ p)&\t{if}\ d=1,\endcases
$$
and hence
$$\sum_{k=1}^{p^a-1}\f{\bar C_k^{(2)}}{m^k}\eq m-6\ (\mo\ p).$$
\endproclaim
\Remark\ 1.2. Let $p>3$ be a prime. By [S98, Corollary 6.1], if $c\in\Z_p$ and $c(c^2+3)\not\eq0\ (\mo\ p)$, then
$c\in C_0(p)\iff u_{(p-(\f p3))/3}\eq0\ (\mo\ p)$,
 where $u_0=0$, $u_1=1$, and $u_{n+1}=6u_n-(3c^2+9)u_{n-1}$ for $n\in\Z^+$.
\medskip

Combining Theorems 1.1-1.3 we obtain the following somewhat surprising result.
\proclaim{Theorem 1.4} Let $p>3$ be a prime. Let $a$ be a positive integer divisible by $6$
and let $d\in\{0,\pm1\}$.
Then $$\sum\Sb 0<k<p^a\\k\eq r\ (\mo\ p-1)\endSb\bi{3k}{k+d}\eq2^{d+3-2r}3^{3r-2}\ (\mo\ p)$$
for all $r\in\Z$, and hence
$$\sum_{k=1}^{p^a-1}\bi{3k}{k+d}\eq-[p=23]3\times2^{d+1}\ (\mo\ p).$$
\endproclaim

We may apply Theorem 1.3 to some particular integers $m=t^2+t+7$ to obtain concrete results.

\proclaim{Theorem 1.5} Let $p\not=3$ be a prime and let $a\in\Z^+$. Then
$$\align\sum_{k=0}^{p^a-1}\f{\bi{3k}k}{9^k}\eq&\cases1&\t{if}\ p^a\eq\pm1\ (\mo\ 9),\\0&\t{if}\ p^a\eq\pm 2\ (\mo\ 9),
\\-1&\t{if}\ p^a\eq\pm4\ (\mo\ 9);\endcases
\\\sum_{k=0}^{p^a-1}\f{\bi{3k}{k-1}}{9^k}\eq&\cases0&\t{if}\ p^a\eq\pm1\ (\mo\ 9),\\1&\t{if}\ p^a\eq\pm 2\ (\mo\ 9),
\\-1&\t{if}\ p^a\eq\pm4\ (\mo\ 9);\endcases
\\\sum_{k=0}^{p^a-1}\f{\bi{3k}{k+1}}{9^k}\eq&\cases0&\t{if}\ p^a\eq\pm1\ (\mo\ 9),\\-5&\t{if}\ p^a\eq\pm 2\ (\mo\ 9),
\\-7&\t{if}\ p^a\eq\pm4\ (\mo\ 9).\endcases
\endalign$$
Consequently,
$$\sum_{k=1}^{p^a-1}\f{C_k^{(2)}}{9^k}\eq-3[p^a\eq\pm2\ (\mo\ 9)]\ (\mo\ p)$$
and
$$\sum_{k=1}^{p^a-1}\f{\bar C_k^{(2)}}{9^k}\eq3[p^a\not\eq\pm1\ (\mo\ 9)]\ (\mo\ p).$$
\endproclaim

\proclaim{Theorem 1.6} Let $p\not=7$ be a prime and let $a\in\Z^+$. Then
$$\align\sum_{k=1}^{p^a-1}\f{\bi{3k}k}{7^k}\eq&-3[p^a\eq\pm2\ (\mo\ 7)]\ (\mo\ p);
\\\sum_{k=0}^{p^a-1}\f{\bi{3k}{k-1}}{7^k}\eq&\cases0&\t{if}\ p^a\eq\pm1\ (\mo\ 7),\\-1&\t{if}\ p^a\eq\pm 2\ (\mo\ 7),
\\1&\t{if}\ p^a\eq\pm3\ (\mo\ 7);\endcases
\\\sum_{k=0}^{p^a-1}\f{\bi{3k}{k+1}}{7^k}\eq&\cases0&\t{if}\ p^a\eq\pm1\ (\mo\ 7),\\-7&\t{if}\ p^a\eq\pm 2\ (\mo\ 7),
\\-1&\t{if}\ p^a\eq\pm3\ (\mo\ 7).\endcases
\endalign$$
Consequently,
$$\sum_{k=0}^{p^a-1}\f{C_k^{(2)}}{7^k}\eq\cases1\ (\mo\ p)&\t{if}\ p^a\eq\pm1\ (\mo\ 7),
\\0\ (\mo\ p)&\t{if}\ p^a\eq\pm2\ (\mo\ 7),
\\-1\ (\mo\ p)&\t{if}\ p^a\eq\pm3\ (\mo\ 7);
\endcases$$
and
$$\sum_{k=1}^{p^a-1}\f{\bar C_k^{(2)}}{7^k}\eq[p^a\not\eq\pm1\ (\mo\ 7)]\ (\mo\ p).$$
\endproclaim

\proclaim{Theorem 1.7} Let $p$ be a prime and let $a\in\Z^+$. If $p\not=5,13$, then
$$\sum_{k=0}^{p^a-1}\f{\bi{3k}k}{13^k}\eq\cases1\ (\mo\ p)&\t{if}\ p^a\eq\pm1,\pm5\ (\mo\ 13),
\\-4/5\ (\mo\ p)&\t{if}\ p^a\eq\pm2,\pm3\ (\mo\ 13),
\\-1/5\ (\mo\ p)&\t{if}\ p^a\eq\pm4,\pm6\ (\mo\ 13),
\endcases$$
and
$$\sum_{k=0}^{p^a-1}\f{\bi{3k}{k+1}}{13^k}\eq\cases1\ (\mo\ p)&\t{if}\ p^a\eq\pm1,\pm5\ (\mo\ 13),
\\-53/5\ (\mo\ p)&\t{if}\ p^a\eq\pm2,\pm3\ (\mo\ 13),
\\-47/5\ (\mo\ p)&\t{if}\ p^a\eq\pm4,\pm6\ (\mo\ 13).
\endcases$$
Also,
$$\sum_{k=0}^{p^a-1}\f{C_k^{(2)}}{13^k}\eq\cases1\ (\mo\ p)&\t{if}\ p^a\eq\pm1,\pm5\ (\mo\ 13),
\\2\ (\mo\ p)&\t{if}\ p^a\eq\pm2,\pm3\ (\mo\ 13),
\\-3\ (\mo\ p)&\t{if}\ p^a\eq\pm4,\pm6\ (\mo\ 13);
\endcases$$
and
$$\sum_{k=0}^{p^a-1}\f{C_k^{(2)}}{19^k}\eq\cases1\ (\mo\ p)&\t{if}\ p^a\eq\pm1,\pm7,\pm8\ (\mo\ 19),
\\-4\ (\mo\ p)&\t{if}\ p^a\eq\pm2,\pm3,\pm5\ (\mo\ 19),
\\3\ (\mo\ p)&\t{if}\ p^a\eq\pm4,\pm6,\pm9\ (\mo\ 19).
\endcases$$
\endproclaim

Now we turn to our results involving third-order and fourth-order Catalan numbers.

\proclaim{Theorem 1.8} Let $p\not=5$ be a prime and let $a\in\Z^+$.
Set
$$S_d=\sum_{k=0}^{p^a-1}\f{\bi{4k}{k+d}}{5^k}\quad \t{for}\ d=-2,-1,\ldots,3p^a.$$

{\rm (i)} When $p\not=11$, we have
$$\align&S_0\eq\cases1\ (\mo\ p)&\t{if}\ p^a\eq1\ (\mo\ 5),\\-9/11\ (\mo\ p)&\t{if}\ p^a\eq -1\ (\mo\ 5),
\\-1/11\ (\mo\ p)&\t{if}\ p^a\eq \pm2\ (\mo\ 5);\endcases
\\&S_1\eq\cases0\ (\mo\ p)&\t{if}\ p^a\eq1\ (\mo\ 5),\\-5/11\ (\mo\ p)&\t{if}\ p^a\eq -1\ (\mo\ 5),
\\-14/11\ (\mo\ p)&\t{if}\ p^a\eq \pm2\ (\mo\ 5);\endcases
\\&S_{-1}\eq\cases 0\ (\mo\ p)&\t{if}\ p^a\eq1\ (\mo\ 5),
\\-3/11\ (\mo\ p)&\t{if}\ p^a\eq -1\ (\mo\ 5),\\7/11\ (\mo\ p)&\t{if}\ p^a\eq 2\ (\mo\ 5),
\\-4/11\ (\mo\ p)&\t{if}\ p^a\eq -2\ (\mo\ 5);\endcases
\\&S_{-2}\eq\cases 0\ (\mo\ p)&\t{if}\ p^a\eq1\ (\mo\ 5),
\\-1/11\ (\mo\ p)&\t{if}\ p^a\eq -1\ (\mo\ 5),\\-16/11\ (\mo\ p)&\t{if}\ p^a\eq 2\ (\mo\ 5),
\\17/11\ (\mo\ p)&\t{if}\ p^a\eq -2\ (\mo\ 5).\endcases
\endalign$$

{\rm (ii)} For $d=2,\ldots,3p^a$ we have
$$S_{d}-S_{d-1}+6S_{d-2}+4S_{d-3}+S_{d-4}\eq\cases6\ (\mo\ p)&\t{if}\ d=p^a+1,
\\4\ (\mo\ p)&\t{if}\ d=2p^a+1,\\0\ (\mo\ p)&\t{otherwise}.\endcases$$

{\rm (iii)} We have
$$\sum_{k=0}^{p^a-1}\f{C_k^{(3)}}{5^k}\eq\cases1\ (\mo\ p)&\t{if}\ p^a\eq1,-2\ (\mo\ 5),
\\0\ (\mo\ p)&\t{if}\ p^a\eq-1\ (\mo\ 5),
\\-2\ (\mo\ p)&\t{if}\ p^a\eq2\ (\mo\ 5).\endcases
$$
Also,
$$\sum_{k=0}^{p^a-1}\f{\bar C_k^{(3)}}{5^k}\eq\cases3\ (\mo\ p)&\t{if}\ p^a\eq1\ (\mo\ 5),
\\-2\ (\mo\ p)&\t{if}\ p^a\eq-1\ (\mo\ 5),
\\1\ (\mo\ p)&\t{if}\ p^a\eq\pm2\ (\mo\ 5).\endcases
$$
\endproclaim

\proclaim{Theorem 1.9} Let $p>3$ be a prime and let $a\in\Z^+$. Then
$$\sum_{k=1}^{p^a-1}\f{3^{3k}}{4^{4k}}C_k^{(3)}\eq\f{(\f{-2}{p^a})-1}{12}\ (\mo\ p)$$
and
$$\sum_{k=1}^{p^a-1}\f{3^{3k}}{4^{4k}}\bi{4k}{k+p^a}\eq-\f{(\f{-2}{p^a})+20}{48}\ (\mo\ p).$$
\endproclaim

\proclaim{Theorem 1.10} Let $p>3$ be a prime.

{\rm (i)} If $(\f p7)=1$, then
$$\sum_{k=1}^{p-1}\f{\bar C_k^{(3)}}{3^k}\eq\cases-6\ (\mo\ p)&\t{if}\ p\eq2\ (\mo\ 3),
\\0\ (\mo\ p)&\t{if}\ p=x^2+3y^2\ \t{and}\ (\f{x+5y}p)=(\f{x-3y}p),
\\-3\ (\mo\ p)&\t{otherwise}.\endcases$$

{\rm (ii)} Suppose that $(\f p{23})=1$.
In the case $p\eq1\ (\mo\ 3)$, if there exists an integer $t\in\Z$ such that $t^2\eq 69\ (\mo\ p)$
and $(97-3t)/2$ is a cubic residue modulo $p$ then
$$\sum_{k=1}^{p-1}(-1)^k\bar C_k^{(4)}\eq0\ (\mo\ p),$$
otherwise
$$\sum_{k=1}^{p-1}(-1)^k\bar C_k^{(4)}\eq-13\ (\mo\ p).$$
In the case $p\eq2\ (\mo\ 3)$, if $v_{(p+1)/3}\eq-13\ (\mo\ p)$
(where $v_0=2,\ v_1=-97$ and $v_{n+1}=-97v_n-13^2v_{n-1}$ for $n\in\Z^+$),
then
$$\sum_{k=1}^{p-1}(-1)^k\bar C_k^{(4)}\eq-10\ (\mo\ p);$$
otherwise we have
$$\sum_{k=1}^{p-1}(-1)^k\bar C_k^{(4)}\eq3\ (\mo\ p).$$
\endproclaim

In the next section we are going to establish a general theorem relating $\sum_{k=0}^{p^a-1}\bi{(h+1)k}{k+d}$ mod $p$
to a linear recurrence of order $h+1$. In Section 3 we shall prove Theorem 1.1.
Theorems 1.2-1.6 will be proved in Section 4. (We omit the proof of Theorem 1.7 since it
is similar to that of Theorem 1.6.) Section 5 is devoted to the proof of Theorem 1.8.
In Section 6 we will show Theorem 1.9. The proof of Theorem 1.10 is very technical, so we omit it.

\heading{2. A general theorem}\endheading

The following lemma is a well known result due to Sylvester which follows from Lagrange's interpolation formula.

\proclaim{Lemma 2.1} Define an $m$-th linear recurrence $\{u_n\}_{n\in\Z}$ by
$$u_0=\cdots=u_{m-2}=0,\ u_{m-1}=1,$$
 and
$$u_{n+m}+a_1u_{n+m-1}+\cdots+a_mu_n=0\quad\ (n\in\Z),$$
where $a_1,\ldots,a_m\in\C$ and $a_m\not=0$. Suppose that
the equation $x^m+a_1x^{m-1}+\cdots+a_0=0$ has
$m$ distinct zeroes $\al_1,\ldots,\al_n\in\C$. Then
$$u_n=\sum_{i=1}^m\f{\al_i^n}{\prod_{j\not=i}(\al_i-\al_j)}\quad\t{for all}\ n\in\Z.$$
\endproclaim

Now we present our general theorem on connections between sums involving binomial coefficients and linear recurrences.
\proclaim{Theorem 2.1}
Let $p$ be a prime and $m\in\Z_p$ with $m\not\eq0\ (\mo\ p)$. Let $a,h\in\Z^+$.
Define an integer sequence $\{u_n\}_{n\in\Z}$ by
$$u_0=\cdots=u_{h-1}=0,\ u_h=1\tag2.1$$
and
$$\sum_{j=0}^{h+1}\(\bi{h+1}j-m\da_{j,h}\)u_{n+j}=0\ \ (n\in\Z).\tag2.2$$

{\rm (i)} For $d\in\{-h+1,\ldots,hp^a\}$ we have
$$\aligned&\sum_{j=0}^{h+1}\(\bi{h+1}j-m\da_{j,h}\)\sum_{k=0}^{p^a-1}\f{\bi{(h+1)k}{k+d+j}}{m^k}
\\\eq&[p^a\mid d+h]\bi{h+1}{(d+h)/p^a+1}\ (\mo\ p)
\endaligned\tag2.3$$
and
$$\sum_{k=0}^{p^a-1}\f{\bi{(h+1)k}{k+d}}{m^k}\eq -\sum_{r=1}^h\bi{h+1}{r+1}u_{h-1+\min\{d-rp^a,0\}}\ (\mo\ p).\tag2.4$$

{\rm (ii)} Suppose that
$$D((1+x)^{h+1}-mx^h)\not\eq0\ (\mo\ p).$$
Then, for $d\in\{-h+1,\ldots,hp^a\}$ we have
$$\aligned\sum_{k=0}^{p^a-1}\f{\bi{(h+1)k}{k+d}}{m^k}\eq&(h+1-m)u_{d+h-1}+u_{p^a+d+h-1}
\\&+\sum_{0<r\ls\lfloor(d-1)/p^a\rfloor}\bi{h+1}{r+1}u_{d+h-1-rp^a}\ (\mo\ p).
\endaligned\tag2.5$$
\endproclaim

\Proof. (i) We first show (2.3) for any given $d\in\{-h+1,\ldots,hp^a\}$.
Observe that
$$\align&\f{\bi{(h+1)p^a}{p^a+d+h}}{m^{p^a-1}}+m\sum_{k=0}^{p^a-1}\f{\bi{(h+1)k}{k+d+h}}{m^k}
\\=&\sum_{k=1}^{p^a}\f{\bi{(h+1)k}{k+d+h}}{m^{k-1}}=\sum_{k=0}^{p^a-1}\f{\bi{(h+1)k+h+1}{k+d+h+1}}{m^k}
\\=&\sum_{k=0}^{p^a-1}\f{\sum_{i=0}^{h+1}\bi{h+1}i\bi{(h+1)k}{k+d+h+1-i}}{m^k}
\\&\quad (\t{by the Chu-Vandermonde identity (see (5.22) of [GKP, p.\,169])})
\\=&\sum_{j=0}^{h+1}\bi{h+1}j\sum_{k=0}^{p^a-1}\f{\bi{(h+1)k}{k+d+j}}{m^k}
\endalign$$
and hence
$$\sum_{j=0}^{h+1}\(\bi{h+1}j-m\da_{j,h}\)\sum_{k=0}^{p^a-1}\f{\bi{(h+1)k}{k+d+j}}{m^k}
\eq\bi{(h+1)p^a}{p^a+d+h}\ (\mo\ p)$$
by Fermat's little theorem.
If $d+h\not\eq0\ (\mo\ p^a)$, then
$$\bi{(h+1)p^a}{p^a+d+h}=\f{(h+1)p^a}{p^a+d+h}\bi{(h+1)p^a-1}{p^a+d+h-1}\eq0\ (\mo\ p);$$
if $d+h=p^aq$ for some $q\in\Z^+$, then
$$\bi{(h+1)p^a}{p^a+d+h}=\bi{(h+1)p^a}{(q+1)p^a}\eq\bi{h+1}{q+1}\ (\mo\ p)$$
by Lucas' theorem (see, e.g., [HS]). Therefore (2.3) follows from the above.

Next we want to prove (2.4) by induction.

For $d\in\{hp^a-h,\ldots,hp^a\}$, as $d\gs h(p^a-1)$ and
$(h-1)p^a-d\ls h-p^a<h$ we have
$$\sum_{k=0}^{p^a-1}\f{\bi{(h+1)k}{k+d}}{m^k}=\f{\bi{(h+1)(p^a-1)}{p^a-1+d}}{m^{p^a-1}}
\eq\da_{d,h(p^a-1)}\ (\mo\ p)$$
and also
$$\align&\sum_{i=1}^h\bi{h+1}{i+1}u_{h-1+\min\{d-ip^a,0\}}
\\=&\sum\Sb 1\ls i\ls h\\ip^a\gs d+h\endSb\bi{h+1}{i+1}u_{h-1+d-ip^a}
\\=&[hp^a\gs d+h]u_{h-1+d-hp^a}=\da_{d,hp^a-h}u_{-1}=-\da_{d,h(p^a-1)}.
\endalign$$
So (2.4) holds for all $d=hp^a-h,\ldots,hp^a$.

Let $-h<d<hp^a-h$ and assume that (2.4) with $d$ replaced by a large integer not exceeding $hp^a$ holds.
For $r\in\{1,\ldots,h\}$, if $ip^a<d+h$ then
$$\sum_{j=0}^{h+1}\(\bi{h+1}j-m\da_{j,h}\)u_{h-1+\min\{d+j-rp^a,0\}}=0$$
since $u_0=\cdots=u_{h-1}=0$;
if $ip^a\gs d+h$, then
$$\align&\sum_{j=0}^{h+1}\(\bi{h+1}j-m\da_{j,h}\)u_{h-1+\min\{d+j-rp^a,0\}}
\\=&\sum_{j=0}^{h}\(\bi{h+1}j-m\da_{j,h}\)u_{h-1+d+j-rp^a}+u_{h-1+\min\{d+h+1-rp^a,0\}}
\\=&\sum_{j=0}^{h}\(\bi{h+1}j-m\da_{j,h}\)u_{h-1+d+j-rp^a}-\da_{d+h,rp^a}=-\da_{d+h,rp^a}.
\endalign$$
So we have
$$\align&\sum_{j=0}^{h+1}\(\bi{h+1}j-m\da_{j,h}\)\sum_{r=1}^h\bi{h+1}{r+1}u_{h-1+\min\{d+j-rp^a,0\}}
\\=&\sum_{i=1}^h\bi{h+1}{r+1}\sum_{j=0}^{h+1}\(\bi{h+1}j-m\da_{j,h}\)u_{h-1+\min\{d+j-rp^a,0\}}
\\=&\sum_{r=1}^h\bi{h+1}{r+1}(-\da_{rp^a,d+h})=-[p^a\mid d+h]\bi{h+1}{(d+h)/p^a+1}.
\endalign$$
Combining this with (2.3) and the induction hypothesis,
we obtain (2.4). This concludes the induction step.

\medskip

(ii)  Write
 $$\sum_{j=0}^{h+1}\(\bi{h+1}j-m\da_{j,h}\)x^j=(x+1)^{h+1}-mx^h=\prod_{i=1}^{h+1}(x-\al_i)$$
 with $\al_1,\ldots,\al_{h+1}\in\C$. As $D:=D((x+1)^{h+1}-mx^h)\not=0$, $\al_1,\ldots,\al_{h+1}$ are distinct.
 Clearly all those $\al_i$, $\al_i^{-1}$, and
 $$c_i:=\f D{\prod_{j\not=i}(\al_i-\al_j)}=\prod\Sb 1\ls s<t\ls h+1\\ s,t\not=i\endSb
 (\al_s-\al_t)^2\times\prod_{j\not=i}(\al_i-\al_j)$$
are  algebraic integers.

 Fix $d\in\{-h+1,\ldots,dp^a\}$. By part (i),
 $$-\sum_{k=0}^{p^a-1}\f{\bi{(h+1)k}{k+d}}{m^k}\eq\sum\Sb 1\ls r\ls h\\ rp^a\gs
  d\endSb\bi{h+1}{r+1}u_{h-1+d-rp^a}\ (\mo\ p).$$
 By Lemma 2.1, for any $n\in\N$ we have
 $$u_n=\sum_{i=1}^{h+1}\f{\al_i^n}{\prod_{j\not=i}(\al_i-\al_j)}
=\f1D\sum_{i=1}^{h+1}c_i\al_i^{n}.$$
Therefore
$$-\sum_{k=0}^{p^a-1}\f{\bi{(h+1)k}{k+d}}{m^k}\eq\f1D\sum_{i=1}^{h+1}c_i\al_i^{d+h-1}
\sum\Sb 1\ls r\ls h\\rp^a\gs d\endSb\bi{h+1}{r+1}\al_i^{-rp^a}
\ (\mo\ p).$$

Since
 $$\sum_{j=0}^{h+1}\bi{h+1}j\al_i^{jp^a}
 \eq\(\sum_{j=0}^{h+1}\bi{h+1}j\al_i^j\)^{p^a}=(m\al_i^h)^{p^a}\eq m\al_i^{hp^a}\ (\mo\ p),$$
 we have
 $$m\eq\sum_{j=0}^{h+1}\bi{h+1}j\al_i^{(j-h)p^a}=\sum_{r=-1}^{h}\bi{h+1}{r+1}\al_i^{-rp^a}\ (\mo\ p)$$
 and hence
 $$\sum_{r=1}^{h}\bi{h+1}{r+1}\al_i^{-rp^a}\eq m-h-1-\al_i^{p^a}\ (\mo\ p).$$
 Therefore
$\sum_{k=0}^{p^a-1}\bi{(h+1)k}{k+d}/m^k$
is congruent to
$$\align&\f1D\sum_{i=1}^{h+1}c_i\al_i^{d+h-1}\(h+1-m+\al_i^{p^a}+\sum_{0<rp^a\ls d-1}\bi{h+1}{r+1}\al_i^{-rp^a}\)
\\=&(h+1-m)u_{d+h-1}+u_{p^a+d+h-1}
+\sum_{0<r\ls\lfloor(d-1)/p^a\rfloor}\bi{h+1}{r+1}u_{d+h-1-rp^a}
\endalign$$
modulo $p$. This proves (2.5).

 The proof of Theorem 2.1 is now complete. \qed

 \heading{3. Proof of Theorem 1.1}\endheading

 To prove Theorem 1.1 in the case $c\eq -1/4\ (\mo\ p)$, we give the following theorem.

 \proclaim{Theorem 3.1} Let $p>3$ be a prime and let $a\in\Z^+$. Then
$$\aligned&\sum_{k=0}^{p^a-1}\f{4^k}{27^k}\bi{3k}{k+d}
\\\eq&\cases ((-1)^d4^{2-d}-7(9d+1)2^d)/81\ (\mo\ p)&\t{if}\ d\in\{-1,\ldots,p^a\},
\\((-1)^d4^{3-d}-(9d+1)2^d)/81\ (\mo\ p)&\t{if}\ d\in\{p^a,\ldots,2p^a\}.\endcases
\endaligned\tag3.1$$
In particular,
$$\align\sum_{k=0}^{p^a-1}\f{4^k}{27^k}\bi{3k}k\eq&\f19\ (\mo\ p),
\quad\ \sum_{k=1}^{p^a-1}\f{4^k}{27^k}\bi{3k}{k+p^a}\eq-\f29\ (\mo\ p),
\\\sum_{k=1}^{p^a-1}\f{4^k}{27^k}\bi{3k}{k+1}\eq&-\f{16}9\ (\mo\ p),
\ \sum_{k=1}^{p^a-1}\f{4^k}{27^k}\bi{3k}{k-1}\eq-\f{4}9\ (\mo\ p).
\endalign$$
\endproclaim
\Proof. Let $u_0=u_1=0,\ u_2=1,$
and
$$u_{n+3}+\l(3-\f{27}4\r)u_{n+1}+u_n=0\quad \ \t{for}\ n=0,1,2,\ldots.$$
 Since
 $$x^3+\l(3-\f{27}4\r)x^2+3x+1=\l(x+\f14\r)(x-2)^2,$$
 there are $a,b,c\in\C$ such that $u_n=(an+b)2^n+c(-1/4)^n$ for all $n\in\N$.
 By $u_0=u_1=0$ and $u_2=1$, we can easily determine the values of $a,b,c$ explicitly. It follows that
 $$u_n=\f{16}{81}\(\l(-\f14\r)^n+\l(\f 98n-1\r)2^n\)\quad\t{for all}\ n\in\N.\tag3.2$$

 Let $d\in\{-1,\ldots,2p^a\}$. Applying (2.4) with $h=2$ and $m=27/4$ we get
 $$\align-\sum_{k=0}^{p^a-1}\f{4^k}{27^k}\bi{3k}{k+d}\eq&\sum_{r=1}^2\bi 3{r+1}u_{1+\min\{d-rp^a,0\}}
 \\\eq&3[d\ls p^a]u_{1+d-p^a}+u_{1+d-2p^a}\ (\mo\ p).
 \endalign$$
 By (3.2) and Fermat's little theorem,
 $$u_{d+1-p^a}\eq\f{(-1)^d4^{2-d}+(9d+1)2^{d+1}}{81}\ (\mo\ p)$$
 and
 $$u_{d+1-2p^a}\eq\f{(-1)^{d-1}4^{3-d}+(9d+1)2^{d}}{81}\ (\mo\ p).$$
 Thus (3.1) follows.

 Applying (3.1) with $d=0,\pm1,p^a$ we immediately obtain the last four congruences in Theorem 3.1.
 We are done. \qed

 \medskip

 Now we need some knowledge about Lucas sequences.

 Given $A,B\in\C$ with $B\not=0$, the Lucas sequences $u_n=u_n(A,B)$ and $v_n=v_n(A,B)$ ($n\in\Z$)
 are defined as follows:
 $$\align &u_0=0,\ u_1=1,\ \t{and}\ u_{n+1}=Au_n-Bu_{n-1}\ (n\in\Z);
\\&v_0=2,\ v_1=A,\ \t{and}\ v_{n+1}=Av_n-Bv_{n-1}\ (n\in\Z).
\endalign$$
It ie easy to see that $v_n=2u_{n+1}-Au_n$ for all $n\in\Z$.
Let $\al$ and $\beta$ be the two roots of the equation $x^2-Ax+B=0$.
It is well known that
$$(\al-\beta)u_n=\al^n-\beta^n\quad\ \t{and}\quad\ v_n=\al^n+\beta^n.$$

\proclaim{Lemma 3.1} Let $p$ be an odd prime and let $a\in\Z^+$. Let $A,B\in\Z_p$ with $\Delta=A^2-4B\not\eq0\ (\mo\ p)$.
Then for any $n\in\Z$ we have
$$u_{n+p^a}\eq\f{Au_n+(\f{\Delta}{p^a})v_n}2\ (\mo\ p)
\ \t{and}\ Bu_{n-p^a}\eq\f{Au_n-(\f{\Delta}{p^a})v_n}2\ (\mo\ p),$$
where $u_k=u_k(A,B)$ and $v_k=v_k(A,B)$.
\endproclaim
\Proof. Let $\al$ and $\beta$ be the two roots of the equation $x^2-Ax+B=0$. Clearly
$$v_{p^a}=\al^{p^a}+\beta^{p^a}\eq(\al+\beta)^{p^a}=A^{p^a}\eq A\ (\mo\ p).$$
Since
$$(\al-\beta)u_{p^a}=\al^{p^a}-\beta^{p^a}\eq(\al-\beta)^{p^a}\ (\mo\ p),$$
we have
$$\Delta u_{p^a}\eq(\al-\beta)^{p^a+1}=\Delta^{(p^a-1)/2}\Delta\ (\mo\ p)$$
and hence
$$u_{p^a}\eq(\Delta^{(p-1)/2})^{\sum_{i=0}^{a-1}p^i}\eq\l(\f{\Delta}{p^a}\r)^a=\l(\f{\Delta}{p^a}\r)\ (\mo\ p).$$

Now,
$$\align 2u_{n+p^a}=&\f{\al^n-\beta^n}{\al-\beta}(\al^{p^a}+\beta^{p^a})+\f{\al^{p^a}-\beta^{p^a}}{\al-\beta}(\al^n+\beta^n)
\\=&u_nv_{p^a}+u_{p^a}v_n\eq Au_n+\l(\f{\Delta}p\r)v_n\ (\mo\ p).
\endalign$$
Also,
$$\align 2u_{n-p^a}=&\f{\al^n-\beta^n}{\al-\beta}(\al^{-p^a}+\beta^{-p^a})
+\f{\al^{-p^a}-\beta^{-p^a}}{\al-\beta}(\al^n+\beta^n)
\\=&u_n\f{\al^{p^a}+\beta^{p^a}}{(\al\beta)^{p^a}}+\f{\beta^{p^a}-\al^{p^a}}{\al-\beta}\cdot\f{v_n}{(\al\beta)^{p^a}}
=u_n\f{v_{p^a}}{B^{p^a}}-\f{u_{p^a}}{B^{p^a}}v_n
\endalign$$
and hence
$$2Bu_{n-p^a}\eq 2B^{p^a}u_{n-p^a}=u_nv_{p^a}-u_{p^a}v_n\eq Au_n-\l(\f{\Delta}{p^a}\r)v_n\ (\mo\ p).$$
This concludes the proof. \qed

\medskip

For Theorem 1.1 in the case $c\not\eq-1/4\ (\mo\ p)$,
we need the following general result.

 \proclaim{Theorem 3.2} Let $p$ be an odd prime and let $a\in\Z^+$.  Let $c\in\Z_p$ with $c\not\eq0,-1,2,-1/4\ (\mo\ p)$,
 and let $d\in\{-1,0,\ldots,p^a\}$. Then
$$\aligned&\sum_{k=0}^{p^a-1}\f{c^{2k}}{(c+1)^{3k}}\bi{3k}{k+d}
\\\eq&u_{d+1}+\f{3c+1}{(c+1)^2(c-2)}\l(u_{d+1}-c^d+\f{u_d}{c^2}\r)
\\&+\f{v_d+c^2v_{d+1}}{2(c+1)^2(c-2)}\(1-\l(\f{4c+1}{p^a}\r)\)\ (\mo\ p),
\endaligned\tag3.3$$
where $u_n=u_n((3c+1)/c^2,-1/c)$ and $v_n=v_n((3c+1)/c^2,-1/c)$.
\endproclaim
\Proof. Set $m=(c+1)^3/c^2$. Then $c$ is a zero of the polynomial
$$x^2+(3-m)x^2+3x+1=(x+1)^3-mx^2.$$
The discriminant of this polynomial is $D=(4m-27)m^2$.
Note that
$$c^2(4m-27)=4(c+1)^3-27c^2=(4c+1)(c-2)^2\not\eq0\ (\mo\ p).$$
We can write
$$x^2+(3-m)x^2+3x+1=(x-c)(x-\al)(x-\beta)$$
with $\al,\beta,c$ distinct. Clearly $-c-\al-\beta=3-m$ and $(-c)(-\al)(-\beta)=1$.
It follows that $\al+\beta=A$ and $\al\beta=B$, where $A=(3c+1)/c^2$ and $B=-1/c$.

Let $U_0=U_1=0$, $U_2=1$ and $U_{n+3}+(3-m)U_{n+2}+3U_{n+1}+U_n=0$ for $n\in\Z$.
Also set $u_n=u_n(A,B)$ and $v_n=v_n(A,B)$ for $n\in\Z$.
By Lemma 2.1, for any $n\in\Z$ we have
$$\align U_n=&\f{c^n}{(c-\al)(c-\beta)}+\f{\al^n}{(\al-c)(\al-\beta)}+\f{\beta^n}{(\beta-c)(\beta-\al)}
\\=&\f1{(c-\al)(c-\beta)}\(c^n+\f{\al^n(\beta-c)-(\al-c)\beta^n}{\al-\beta}\)
\\=&\f1{c^2-Ac+B}\(c^n+B\f{\al^{n-1}-\beta^{n-1}}{\al-\beta}-c\f{\al^n-\beta^n}{\al-\beta}\)
\\=&\f c{c^3-3c-2}(c^n-c^{-1}u_{n-1}-cu_n)=\f{c^{n+1}-u_{n-1}-c^2u_n}{(c+1)^2(c-2)}.
\endalign$$

In light of Theorem 2.1(ii),
$$\sum_{k=0}^{p^a-1}\f{\bi{3k}{k+d}}{m^k}
\eq(3-m)U_{d+1}+U_{p^a+d+1}\ (\mo\ p)$$
and hence
$$\aligned&\sum_{k=0}^{p^a-1}\f{c^{2k}}{(c+1)^{3k}}\bi{3k}{k+d}-U_{p^a+d+1}
\\\eq&\l(3-\f{(c+1)^3}{c^2}\r)\f{c^{d+2}-u_d-c^2u_{d+1}}{(c+1)^2(c-2)}\ (\mo\ p).
\endaligned\tag3.4$$

 Note that
$$\Delta:=A^2-4B=\f{(3c+1)^2}{c^4}+\f4c=\f{(c+1)^2(4c+1)}{c^4}\not\eq0\ (\mo\ p)$$
and
$$\l(\f{\Delta}{p^a}\r)=\l(\f{4c+1}{p^a}\r).$$
By Lemma 3.1,
$$2u_{p^a+d}\eq\l(\f{4c+1}p\r)v_d+Au_d\ (\mo\ p)$$
and
$$2u_{p^a+d+1}\eq\l(\f{4c+1}p\r)v_{d+1}+Au_{d+1}\ (\mo\ p).$$
Thus
$$\align &U_{p^a+d+1}=\f{c^{p^a+d+2}-u_{p^a+d}-c^2u_{p^a+d+1}}{(c+1)^2(c-2)}
\\\eq&\f{2c^{d+3}-((\f{4c+1}p)v_d+Au_d)-c^2((\f{4c+1}p)v_{d+1}+Au_{d+1})}{2(c+1)^2(c-2)}\ (\mo\ p).
\endalign$$
Note that
$$\align&\f{v_d+c^2v_{d+1}}2+A\f{u_d+c^2u_{d+1}}2
\\=&\f{v_d+Au_d}2+c^2\f{v_{d+1}+Au_{d+1}}2=u_{d+1}+c^u_{d+2}
\\=&u_{d+1}+c^2\l(\f{3c+1}{c^2}u_{d+1}+\f{u_d}c\r)=(3c+2)u_{d+1}+cu_d.
\endalign$$
Therefore
$$U_{p^a+d+1}\eq\f{c^{d+3}+(v_d+c^2v_{d+1})\f{1-(\f{4c+1}{p^a})}2-((3c+2)u_{d+1}+cu_d)}{(c+1)^2(c-2)}\ (\mo\ p).$$
Combining this with (3.4) we finally obtain the desired (3.3). \qed

\proclaim{Corollary 3.1} Let $p>3$ be a prime and let $d\in\{-1,0,\ldots,p^a\}$ with $a\in\Z^+$. Then
$$\sum_{k=0}^{p^a-1}\f{3^k}{8^k}\bi{3k}{k+d}
\eq\cases\f{(-3)^{d/2}}{28}(1+27(\f {p^a}3))\ (\mo\ p),&\t{if}\ 2\mid d,
\\\f{(-3)^{(d+3)/2}}{28}(1-(\f {p^a}3))\ (\mo\ p),&\t{if}\ 2\nmid d.
\endcases\tag3.5$$
\endproclaim
\Proof. Set $c=-1/3$. Then $c^2/(c+1)^3=3/8$, $(3c+1)/c^2=0$ and $-1/c=3$.
Let $u_n=u_n(0,3)$ and $v_n=v_n(0,3)$ for $n\in\Z$.
We clearly have
$$u_{2n}=v_{2n+1}=0,\ u_{2n+1}=(-3)^n\ \t{and}\ v_{2n}=2(-3)^n\ \ \t{for all}\ n\in\Z.$$
Applying Theorem 3.2 we immediately get the desired result. \qed

\medskip
\noindent{\it Proof of Theorem 1.1}. In the case $c\eq-1/4\ (\mo\ p)$, we have $c^2/(c+1)^3\eq 4/27\ (\mo\ p)$
and $c'\eq -8/9\ (\mo\ p)$, hence the desired congruences follow from Theorem 3.1.

Below we assume that $c\not\eq-1/4\ (\mo\ p)$.
For the first three congruences in Theorem 1.1,
we may simply apply Theorem 3.2 with $d=0,\pm1$.

As in the proof of Theorem 3.2, we define $A=(3c+1)/c^2$, $B=-1/c$ and $\Delta=A^2-4B$.
Let $u_n=u_n(A,B)$ and $v_n=v_n(A,B)$ for $n\in\Z$. By Lemma 3.1,
$$2u_{p^a+1}=Au_1+\l(\f{\Delta}{p^a}\r)v_1=A+A\l(\f{4c+1}{p^a}\r)\ (\mo\ p)$$
and
$$\align v_{p^a+1}=&2u_{p^a+2}-Au_{p^a+1}=Au_{p^a+1}-2Bu_{p^a}
\\\eq&\f{A^2+A^2(\f{4c+1}{p^a})}2-2B\l(\f{4c+1}{p^a}\r)
=\f{A^2+\Delta(\f{4c+1}{p^a})}2\ (\mo\ p).
\endalign$$
These, together with Theorem 3.2 in the case $d=p^a$, yield the last congruence in Theorem 1.1.
We are done. \qed

\heading{4. Proofs of Theorems 1.2-1.6}\endheading

\proclaim{Lemma 4.1} Let $p>3$ be a prime and let $a\in\Z^+$. Let
$$u_0=u_1=0,\ u_2=1,\ \t{and}\ u_{n+3}+a_1u_{n+2}+a_2u_{n+1}+a_3u_n=0\ \t{for all}\ n\in\N,$$
where $a_1,a_2,a_3\in\Z$. Suppose that $d\in\Z$ and
$$d^2\eq D(x^3+a_1x^2+a_2x+a_3)\not\eq0\ (\mo\ p).$$ Set
$b=-2a_1^3+9a_1a_2-27a_3$. Then
$$u_{p^a}\eq\cases0\ (\mo\ p)&\t{if}\ p\mid a_1^2-3a_2\ \t{or}\ b/(3d)\in C_0(p^a),
\\\pm(a_1^2-3a_2)/d\ (\mo\ p)&\t{if}\ \pm b/(3d)\in C_1(p^a);\endcases$$
$$u_{p^a+1}\eq\cases b^{(p^a-1)/3}\ (\mo\ p)&\t{if}\ p\mid a_1^2-3a_2,
\\1\ (\mo\ p)&\t{if}\ b/(3d)\in C_0(p^a),
\\(\pm(9a_3-a_1a_2)-d)/(2d)\ (\mo\ p)&\t{if}\ \pm b/(3d)\in C_1(p^a);\endcases$$
$$u_{p^a+2}\eq\cases -a_1(2b^{(p^a-1)/3}+1)/3\ (\mo\ p)&\t{if}\ p\mid a_1^2-3a_2,
\\-a_1\ (\mo\ p)&\t{if}\ b/(3d)\in C_0(p^a),
\\\pm(a_2^2-3a_1a_3)/d\ (\mo\ p)&\t{if}\ \pm b/(3d)\in C_1(p^a).\endcases$$
\endproclaim
\Proof. In the case $a=1$, this is a result due to Z. H. Sun [S03, Theorems 3.2-3.3].
Modifying the proof for the case $a=1$ slightly, we get the result with general $a$. \qed

Actually we just need the following particular result implied by Lemma 4.1.

\proclaim{Lemma 4.2} Let $p>3$ be a prime and let $a\in\Z^+$. Let $m,t\in\Z$ with $2t+1\not\eq0\ (\mo\ p)$
$m\eq t^2+t+7\not\eq0\ (\mo\ p)$.
 Define $\{u_n\}_{n\gs0}$ by
$$u_0=u_1=0,\ u_2=1,\ \t{and}\ u_{n+3}+(3-m)u_{n+2}+3u_{n+1}+u_n=0\ \t{for}\ n\in\N.$$
Set $c=(2m^2-18m+27)/(6t+3)$. Then
$$\align u_{p^a}\eq&\cases0\ (\mo\ p)&\t{if}\ p\mid m-6\ \t{or}\ c\in C_0(p^a),
\\\pm(m-6)/(2t+1)\ (\mo\ p)&\t{if}\ \pm c\in C_1(p^a);\endcases
\\u_{p^a+1}\eq&\cases 2^{(p^a-1)/3}\ (\mo\ p)&\t{if}\ p\mid m-6,
\\1\ (\mo\ p)&\t{if}\ c\in C_0(p^a),
\\\pm 3/(4t+2)-1/2\ (\mo\ p)&\t{if}\ \pm c\in C_1(p^a);\endcases
\\u_{p^a+2}\eq&\cases 2^{(p^a+2)/3}+1\ (\mo\ p)&\t{if}\ p\mid m-6,
\\m-3\ (\mo\ p)&\t{if}\ c\in C_0(p^a),
\\\pm3/(2t+1)\ (\mo\ p)&\t{if}\ \pm c\in C_1(p^a).\endcases
\endalign$$
\endproclaim

\noindent{\it Proof of Theorem 1.2}. The discriminant of the polynomial
$(x+1)^3-6x^2$ is $D=(4\times6-27)6^2=-108$.

{\it Case}\ 1. $p\eq-1\ (\mo\ 3)$. In this case, $(\f Dp)=-1$ and hence
$(x+1)^2-6x^2$ mod $p$ has exactly two irreducible factors, thus
 $(c+1)^3\eq 6c^2\ (\mo\ p)$ for some $c\in\Z$.
Clearly $c\not\eq0,-1,2,-1/4\ (\mo\ p)$. Note that $a$ is even since $p^a\eq1\ (\mo\ 3)$.
As
$$\l(\f{4c+1}{p^a}\r)=\l(\f{4c+1}{p}\r)^a=1,$$
the first congruence in Theorem 1.2 follows from Theorem 1.1.

{\it Case}\ 2. $p\eq 1\ (\mo\ 3)$. In this case, for some $t\in\Z$ we have $(2t+1)^2\eq-3\ (\mo\ p)$, i.e.,
$t^2+t+7\eq 6\ (\mo\ p)$. Let $u_0=u_1=0$, $u_2=1$ and
$$u_{n+3}+(3-6)u_{n+2}+3u_{n+1}+u_n=0\ (n=0,1,2,\ldots).$$
By Theorem 2.1(ii), for $d=-1,\ldots,p^a$ we have
$$\sum_{k=0}^{p^a-1}\f{\bi{3k}{k+d}}{6^k}\eq u_{p^a+d+1}+(3-6)u_{d+1}\ (\mo\ p).$$
Combining this with Lemma 4.2 in the case $m=6$, we are able to determine $\sum_{k=0}^{p^a-1}\bi{3k}{k+d}/6^k$ mod $p$
for $d=0,\pm1$. Note that
$$\bar C_k^{(2)}=\f2{k+1}\bi{3k}k=2\bi{3k}k-\bi{3k}{k+1}.$$
So we have all the desired congruences in Theorem 1.2. \qed

\medskip
\noindent{\it Proof of Theorem 1.3}. Define $\{u_n\}_{n\gs0}$ as in Lemma 4.2. By Theorem 2.1(ii),
for $d=-1,\ldots,p^a$ we have
$$\sum_{k=0}^{p^a-1}\f{\bi{3k}{k+d}}{m^k}\eq u_{p^a+d+1}+(3-m)u_{d+1}\ (\mo\ p).$$
Observe that
$c^2+3\not\eq0\ (\mo\ p)$ since
$$(2m^2-18m+27)^2+3(6t+3)^2\eq4m(m-6)^2\ (\mo\ p).$$
By applying Lemma 4.2 we obtain the desired result. \qed

\medskip
\noindent{\it Proof of Theorem 1.4}. Fix $d\in\{0,\pm1\}$.
By Theorem 3.1 (or Theorem 1.1 in the case $c=-1/4$),
$$\sum_{k=1}^{p^a-1}\f{4^k}{27^k}\bi{3k}{k+d}\eq-\f{2^{d+3}}9\ (\mo\ p).$$

Let $m\in\{1,\ldots,p-1\}$ with $m\not\eq 27/4\ (\mo\ p)$.
If $(c+1)^3\eq m c^2\ (\mo\ p)$ for some $c\in\Z$,
then $c\not\eq0,-1,2,-1/4\ (\mo\ p)$. Thus, by Theorem 1.1 we have
$$\sum_{k=1}^{p^a-1}\f{\bi{3k}{k+d}}{m^k}\eq0\ (\mo\ p)$$
since $(\f{4c+1}{p^a})=(\f{4c+1}p)^a=1$.

Now assume that $(x+1)^3\eq 6 x^2\ (\mo\ p)$ is not solvable over $\Z$. Then, by Stickelberger's theorem,
$$\l(\f p3\r)=\l(\f{-108}p\r)=\l(\f{D((1+x)^3-3x^2)}p\r)\eq(-1)^3-1=1$$
and hence $p\eq1\ (\mo\ 3)$. By Theorem 1.2,
$$\sum_{k=1}^{p^a-1}\f{\bi{3k}{k+d}}{6^k}\eq0\ (\mo\ p)$$
since
$$2^{(p^6-1)/3}=2^{(p^2-1)(p^4+p^2+1)/3}\eq1\ (\mo\ p).$$

Now suppose that $m\not\eq6\ (\mo\ p)$ and $(x+1)^3\eq m x^2\ (\mo\ p)$ is not solvable over $\Z$.
Then
$$\l(\f{(4m-27)m^2}p\r)=\l(\f{D((x+1)^3-mx^2)}p\r)=1$$
and hence $m\eq t^2+t+7\ (\mo\ p)$ for some $t\in\Z$ with $t\not\eq-1/2\ (\mo\ p)$.
Let $c=(2m^2-18m+27)/(6t+3)$. By Theorem 1.3,
$$\l(\f{c+1+2\omega}{p^a}\r)_3=\l(\f{c+1+2\omega}{p}\r)_3^3=1.$$
Hence $c\in C_0(p^a)$ and $$\sum_{k=1}^{p^a-1}\f{\bi{3k}{k+d}}{m^k}\eq0\ (\mo\ p).$$

In view of the above,
$$\align-\sum\Sb 0<k<p^a\\k\eq r\ (\mo\ p-1)\endSb\bi{3k}{k+d}\eq
&\sum_{k=1}^{p^a-1}\bi{3k}{k+d}(p-1)[p-1\mid k-r]
\\\eq&\sum_{k=1}^{p^a-1}\bi{3k}{k+d}\sum_{m=1}^{p-1}m^{r-k}
=\sum_{m=1}^{p-1}m^r\sum_{k=1}^{p^a-1}\f{\bi{3k}{k+d}}{m^k}
\\\eq&\f{27^r}{4^r}\sum_{k=1}^{p^a-1}\f{4^k}{27^k}\bi{3k}{k+d}\eq -\f{27^r}{4^r}\cdot\f{2^{d+3}}9\ (\mo\ p).
\endalign$$
So we have the first congruence in Theorem 1.4. The second congruence follows immediately
since
$$\sum_{k=1}^{p^a-1}\bi{3k}{k+d}=\sum_{r=0}^{p-2}\sum\Sb 0<k<p^a\\ k\eq r\ (\mo\ p-1)\endSb\bi{3k}{k+d}$$
and
$$\sum_{r=0}^{p-2}\f{27^r}{4^r}=\f{27^{p-1}/4^{p-1}-1}{27/4-1}\eq-[p=23]\ (\mo\ p).$$
This concludes the proof of Theorem 1.4. \qed

\medskip
\noindent{\it Proof of Theorem 1.5}.
It suffices to deduce the first, the second  and the third congruences in Theorem 1.5.
Since we can handle the case $p=2$ by detailed analysis, below we assume $p>3$.

By Theorem 1.3 in the case $m=9$ and $t=1$, we only need to show that
$$\aligned3\in C_0(p^a)&\iff p^a\eq\pm1\ (\mo\ 9),
\\3\in C_1(p^a)&\iff p^a\eq\pm2\ (\mo\ 9),
\\3\in C_2(p^a)&\iff p^a\eq\pm 4\ (\mo\ 9).
\endaligned\tag4.1$$
Note that
$$\l(\f{3+1+2\omega}{p^a}\r)_3=\l(\f2p\r)_3^a\l(\f{2+\omega}{p^a}\r)_3=\l(\f{2+\omega}{p^a}\r)$$
and
$$\overline{\l(\f{2+\omega}{p^a}\r)_3}=\l(\f{2+\bar\omega}{p^a}\r)_3=\l(\f{1-\omega}{p^a}\r)_3
=\omega^{((\f{p^a}3)p^a-1)/3}.$$
(See, e.g., [IR].) Clearly,
$$\f{(\f{p^a}3)p^a-1}3\eq\cases0&\t{if}\ p^a\eq\pm1\ (\mo\ 9),
\\2\ (\mo\ 3)&\t{if}\ p^a\eq\pm2\ (\mo\ 9),\\1\ (\mo\ 3)&\t{if}\ p^a\eq\pm4\ (\mo\ 9).\endcases$$
Therefore the three formulae in (4.1) are valid. We are done. \qed

\medskip
\noindent{\it Proof of Theorem 1.6}. We only need to deduce the first, the second
and the third congruences in Theorem 1.6.
Since we can handle the case $p=2,3$ by detailed analysis, below we assume $p>3$.

By Theorem 1.3 in the case $m=7$ and $t=0$, it suffices to show that
$$\aligned-\f13\in C_0(p^a)&\iff p^a\eq\pm1\ (\mo\ 7),
\\-\f13\in C_1(p^a)&\iff p^a\eq\pm3\ (\mo\ 7),
\\-\f13\in C_2(p^a)&\iff p^a\eq\pm 2\ (\mo\ 7).
\endaligned\tag4.2$$
Clearly
$$\l(\f 3{p^a}\r)_3\l(\f{-1/3+1+2\omega}{p^a}\r)_3=\l(\f2{p^a}\r)_3\l(\f{1+3\omega}{p^a}\r)_3,$$
and hence
$$\l(\f{-1/3+1+2\omega}{p^a}\r)_3=\l(\f{1+3\omega}{p^a}\r)_3$$
since $(\f2{p^a})=(\f3{p^a})=1$.
Observe that the norm of $1+3\omega$ is $N(1+3\omega)=(1+3\omega)(1+3\bar\omega)=7$.
By the cubic reciprocity law,
$$\l(\f{1+3\omega}{p^a}\r)_3=\l(\f {p^a}{1+3\omega}\r)_3.$$
If $p^a\eq\pm1\ (\mo\ 7)$, then
$$\l(\f {p^a}{1+3\omega}\r)_3=\l(\f{\pm1}{1+3\omega}\r)_3=\l(\f{\pm1}{1+3\omega}\r)_3^3=1$$
and hence $-1/3\in C_0(p^a)$.
If $p^a\eq\pm2\ (\mo\ 7)$, then
$$\l(\f {p^a}{1+3\omega}\r)_3=\l(\f{\pm2}{1+3\omega}\r)_3\eq(\pm2)^{(N(1+3\omega)-1)/3}=4\eq\omega^2
\ (\mo\ 1+3\omega),$$
hence $(\f{p^a}{1+3\omega})_3=\omega^2$ and $-1/3\in C_2(p^a)$.
If $p^a\eq\pm4\ (\mo\ 7)$, then
$$\l(\f {p^a}{1+3\omega}\r)_3=\l(\f{\pm4}{1+3\omega}\r)_3=\l(\f2{1+3\omega}\r)_3^2=(\omega^2)^2=\omega$$
and hence $-1/3\in C_1(p^a)$. This completes the proof. \qed

\medskip

 \heading{5. Proof of Theorem 1.8}\endheading

 In this section we define a sequence $\{u_n\}_{n\in\Z}$ by
 $$u_0=u_1=u_2=0,\ u_3=1$$
 and
 $$u_{n+4}-u_{n+3}+6u_{n+2}+4u_{n+1}+u_n=0\ \ (n\in\Z).$$
 We also set
 $$v_n^{(1)}=u_{n+2}-3u_{n+1}\quad\t{and}\quad v_n^{(2)}=3u_{n+1}+2u_n.\tag5.1$$

 Recall that the Lucas sequence $\{L_n\}_{n\in\Z}$ is given by
 $$L_0=2,\ L_1=1, \t{and}\ L_{n+1}=L_n+L_{n-1}\quad\t{for all}\ n\in\Z.$$

 \proclaim{Lemma 5.1} {\rm (i)} We have
 $$x^4-x^3+6x^2+4x+1=(x+1)^4-5x^3=\prod\Sb \zeta^5=1\\\zeta\not=1\endSb(x-(1+\zeta)^2).\tag5.2$$

 {\rm (ii)} Let $p$ be a prime, and let $a\in\Z^+$ and $s\in\{1,2\}$.
 Then, for any $d\in\N$ we have
 $$\align5(v_{p^a+d}^{(s)}-v_d^{(s)})
 \eq&2L_{2d}([5\mid d+2p^a-2s+1]-[5\mid d+2p^a-2s])
 \\&+4L_{2d}([5\mid d+p^a-2s+1]-[5\mid d+p^a-2s])
 \\&+\(\f{d+2p^a-2s+1}5\)L_{2d-(\f{d+2p^a-2s+1}5)}
 \\&-\(\f{d+2p^a-2s}5\)L_{2d-(\f{d+2p^a-2s}5)}
 \\&+2\(\f{d+p^a-2s+1}5\)L_{2d-(\f{d+p^a-2s+1}5)}
 \\&-2\(\f{d+p^a-2s}5\)L_{2d-(\f{d+p^a-2s}5)}\ (\mo\ p).
 \endalign$$
 \endproclaim
 \Proof. (i) It is easy to verify that
 $$(1+(1+x)^2)^4-5(1+x)^6=\f{x^5-1}{x-1}(x^4+7x^3+19x^2+23x+11).$$
 Therefore any primitive 5th root $\zeta$ of unity is a zero of $(1+x)^4=5x^3.$
 So (5.2) follows.

\medskip

(ii) For $n\in\Z$ let
$$V_n^{(s)}=\f15\sum\Sb \zeta^5=1\\\zeta\not=1\endSb(\zeta^{1-2s}-\zeta^{-2s})(1+\zeta)^{2n}
=\f15\sum_{\zeta^5=1}(\zeta^{1-2s}-\zeta^{-2s})(1+\zeta)^{2n}.$$
Then $\{V_n\}_{n\in\Z}$ satisfies the recurrence relation
$$V^{(s)}_{n+4}-V^{(s)}_{n+3}+6V^{(s)}_{n+2}+4V^{(s)}_{n+1}+V^{(s)}_n=0\ \ (n\in\Z).$$
Clearly we also have
$$v^{(s)}_{n+4}-v^{(s)}_{n+3}+6v^{(s)}_{n+2}+4v^{(s)}_{n+1}+v^{(s)}_n=0\ \ (n\in\Z).$$

Note that
$$\f15\sum_{\zeta^5=1}\zeta^k=[5\mid k]\ \quad\t{for any}\ k\in\Z;$$
in particular
$$\f15\sum_{\zeta^5=1}\zeta^{1-2s}=0=\f15\sum_{\zeta^5=1}\zeta^{-2s}.$$
Thus
$$V^{(s)}_0=\f15\(\sum_{\zeta^5=1}\zeta^{1-2s}-\sum_{\zeta^5=1}\zeta^{-2s}\)=0=v_0^{(s)}$$
and
$$\align V^{(s)}_1=&\f15\sum_{\zeta^5=1}(\zeta^{1-2s}-\zeta^{-2s})(1+2\zeta+\zeta^2)
\\=&\f15\sum_{\zeta^5=1}(\zeta^{3-2s}+\zeta^{2-2s})=[s=1]=v_1^{(s)}.
\endalign$$
Also,
$$\align V^{(s)}_2=&\f15\sum_{\zeta^5=1}(\zeta^{1-2s}-\zeta^{-2s})(1+4\zeta+6\zeta^2+4\zeta^3+\zeta^4)
\\=&\f15\sum_{\zeta^5=1}(3\zeta^{4-2s}+2\zeta^{3-2s}-2\zeta^{2-2s})=-2[s=1]+3[s=2]=v_2^{(s)}
\endalign$$
and
$$\align V^{(s)}_3=&\f15\sum_{\zeta^5=1}(\zeta^{1-2s}-\zeta^{-2s})
(1+6\zeta+15\zeta^2+20\zeta^3+15\zeta^4+6\zeta^5+\zeta^6)
\\=&\f15\sum_{\zeta^5=1}(\zeta^{1-2s}-\zeta^{-2s})(7\zeta+15\zeta^2+20\zeta^3+15\zeta^4)
\\=&[s=1](7-15)+[s=2](20-15)=v_3^{(s)}.
\endalign$$

By the above, $V_n^{(s)}=v_n^{(s)}$ for all $n\in\N$.

Now fix $d\in\N$. For any algebraic integer $\zeta$, we have $(1+\zeta)^{p^a}\eq1+\zeta^{p^a}\ (\mo\ p)$
and hence
$$\align&(1+\zeta)^{2(p^a+d)}-(1+\zeta)^{2d}
\\\eq&(1+\zeta)^{2d}((1+\zeta^{p^a})^2-1)
\\\eq&\sum_{k=0}^{2d}\bi{2d}k(\zeta^{k+2p^a}+2\zeta^{k+p^a})\ (\mo\ p).
\endalign$$
Thus
$$\align &5(V^{(s)}_{p^a+d}-V^{(s)}_d)
\\=&\sum_{\zeta^5=1}(\zeta^{1-2s}-\zeta^{-2s})((1+\zeta)^{2p^a+2d}-(1+\zeta)^{2d})
\\\eq&\sum_{\zeta^5=1}(\zeta^{1-2s}-\zeta^{-2s})\sum_{k=0}^{2d}\bi{2d}k(\zeta^{k+2p^a}+2\zeta^{k+p^a})
\\\eq&5\sum_{k+2p^a\eq 2s-1\ (\mo\ 5)}\bi{2d}k-5\sum_{k+2p^a\eq2s\ (\mo\ 5)}\bi{2d}k
\\&+10\sum_{k+p^a\eq 2s-1\ (\mo\ 5)}\bi{2d}k-10\sum_{k+p^a\eq2s\ (\mo\ 5)}\bi{2d}k\ (\mo\ p).
\endalign$$
It is known that
$$5\sum_{k\eq r\ (\mo\ 5)}\bi{2d}k-2^{2d}=[5\mid d-r]2L_{2d}+\l(\f{d-r}5\r)L_{2d-(\f{d-r}5)}$$
for all $r\in\Z$. (Cf. [S92], [SS], [Su02] and [Su08].) Therefore $5(V^{(s)}_{p^a+d}-V^{(s)}_d)$ is congruent to
the right-hand side of the congruence in Lemma 5.1(ii)  modulo $p$. So the desired congruence follows.

 The proof of Lemma 5.1 is now complete.  \qed

 \Remark\ 5.1. On April 27, 2009, the author sent a message [Su09] to Number Theory List
 in which he raised the following conjecture:
Let $p$ be a prime and $N_p$ denote the number of solutions
of the the congruence $x^4-x^3+6x^2+4x+1\eq0\ (\mo\ p)$.
If $p\eq1\ (\mo\ 10)$ and $p\not=11$, then $N_p=4$; if $p\eq3,7,9\ (\mo\ 10)$ then $N_p=0$.
Also,
$$v_p^{(1)}=u_{p+2}-3u_{p+1}\eq\cases 1\ (\mo\ p)&\t{if}\ p\eq1,3\ (\mo\ 10),
\\-2\ (\mo\ p)&\t{if}\ p\eq7\ (\mo\ 10),\\0\ (\mo\ p)&\t{if}\ p\eq9\ (\mo\ 10).
\endcases$$
In May 2009, the conjecture was confirmed by K. Buzzard [B], R. Chapman [Ch], E.H. Goins [G]
and also D. Brink, K. S. Chua, K. Foster and F. Lemmermeyer (personal communications);
all of them realized Lemma 5.1(i).
The author would like to thank these cleaver mathematicians for their solutions to the problem.

 \proclaim{Lemma 5.2} Let $p\not=5$ be a prime and let $a\in\Z^+$.
For $s=1,2$ we have
$$v_{p^a}^{(s)}\eq[5\mid p^a-s-2]-[5\mid p^a-s]+2[5\mid p^a-2s+1]-2[5\mid p^a-2s]\ (\mo\ p).$$
Also,
$$v_{p^a+1}^{(1)}-1\eq\cases-3\ (\mo\ p)&\t{if}\ p^a\eq 1\ (\mo\ 5),
\\2\ (\mo\ p)&\t{if}\ p^a\eq -1\ (\mo\ 5),\\-1\ (\mo\ p)&\t{if}\ p^a\eq\pm2\ (\mo\ 5);\endcases$$
$$v_{p^a+1}^{(2)}\eq\cases\pm3\ (\mo\ p)&\t{if}\ p^a\eq \pm1\ (\mo\ 5),\\\pm1\ (\mo\ p)&\t{if}\ p^a\eq\pm2\ (\mo\ 5);
\endcases$$
$$v_{p^a+2}^{(1)}-v_{2}^{(1)}\eq\cases-6\ (\mo\ p)&\t{if}\ p^a\eq1\ (\mo\ 5),
 \\7\ (\mo\ p)&\t{if}\ p^a\eq -1\ (\mo\ 5),\\2\ (\mo\ p)&\t{if}\ p^a\eq 2\ (\mo\ 5),
 \\3\ (\mo\ p)&\t{if}\ p^a\eq -2\ (\mo\ 5);\endcases$$
$$v_{p^a+2}^{(2)}-v_{2}^{(2)}\eq\cases2\ (\mo\ p)&\t{if}\ p^a\eq1\ (\mo\ 5),
 \\-3\ (\mo\ p)&\t{if}\ p^a\eq -1\ (\mo\ 5),\\-4\ (\mo\ p)&\t{if}\ p^a\eq \pm2\ (\mo\ 5);\endcases$$
 $$v_{p^a+3}^{(2)}-v_{3}^{(2)}\eq\cases-18\ (\mo\ p)&\t{if}\ p^a\eq1\ (\mo\ 5),
 \\16\ (\mo\ p)&\t{if}\ p^a\eq -1\ (\mo\ 5),\\-8\ (\mo\ p)&\t{if}\ p^a\eq 2\ (\mo\ 5),
 \\-5\ (\mo\ p)&\t{if}\ p^a\eq -2\ (\mo\ 5);\endcases$$
and
 $$v_{p^a-1}^{(1)}\eq\cases0\ (\mo\ p)&\t{if}\ p^a\eq\pm1\ (\mo\ 5),
 \\\pm5\ (\mo\ p)&\t{if}\ p^a\eq \pm2\ (\mo\ 5).\endcases$$
\endproclaim
\Proof. Note that for $a\in\Z$ we have
$$\l(\f a5\r)L_{-(\f a5)}=-\l(\f a5\r)^2=-[5\nmid a]=[5\mid a]-1.$$ Thus
Lemma 5.1 in the case $d=0$ yields the first congruence in Lemma 5.2.
We can also apply Lemma with $d=1,2,3$ to get the five congruences in Lemma 5.2 following the first one.

 Now we deduce the last congruence in Lemma 5.2. By the proof of Lemma 5.1,
 $$5v_{p^a-1}^{(1)}=5V_{p^a-1}^{(1)}=\sum_{\zeta^5=1}(\zeta^{-1}-\zeta^{-2})((1+\zeta)^{2(p^a-1)}\ (\mo\ p).$$
For any primitive 5th root $\zeta$ of unity, clearly
$$(1+\zeta)(\zeta+\zeta^3)=\zeta+\zeta^3+\zeta^2+\zeta^4=-1$$
and hence
$$(1+\zeta)^{-2}=(-\zeta-\zeta^3)^2=2\zeta^4+\zeta^2+\zeta=\zeta^4-\zeta^3-1;$$
also
$$(\zeta^{-1}-\zeta^{-2})(\zeta^4-\zeta^3-1)=\zeta-\zeta^{-1}-2\zeta^2+2\zeta^{-2}$$
and
$$(1+\zeta)^{2p^a}\eq(1+\zeta^{p^a})^2\eq1+2\zeta^{p^a}+\zeta^{2p^a}\ (\mo\ p).$$
Therefore
$$\align5v_{p^a-1}^{(1)}\eq&\sum_{\zeta^5=1}(\zeta-\zeta^{-1}-2\zeta^2+2\zeta^{-2})(1+2\zeta^{p^a}+\zeta^{2p^a})
\\\eq&\sum_{\zeta^5=1}(\zeta-\zeta^{-1}-2\zeta^2+2\zeta^{-2})(2\zeta^{p^a}+\zeta^{2p^a})
\\\eq&\cases5((-1)\times2+2\times1)\ (\mo\ p)&\t{if}\ p^a\eq1\ (\mo\ 5),
\\5(1\times2+(-2)\times1)\ (\mo\ p)&\t{if}\ p^a\eq-1\ (\mo\ 5),
\\5(2\times2+1\times 1)\ (\mo\ p)&\t{if}\ p^a\eq 2\ (\mo\ 5),
\\5(-2\times2+(-1)\times 1)\ (\mo\ p)&\t{if}\ p^a\eq -2\ (\mo\ 5).
\endcases\endalign$$
This yields the last congruence in Lemma 5.2. We are done. \qed

\medskip
\noindent{\it Proof of Theorem 1.8}. For the polynomial
$$x^4-x^3+6x^2+4x+1=(x+1)^4-5x^3,$$
its discriminant is $5^3\times 11^2$.

(i) Suppose that $p\not=11$. Then $p$ does not divide $D((x+1)^4-5x^3)$.
For any $n\in\Z$ we have
$$11u_n=(3u_{n+1}+2u_n)-3(u_{n+1}-3u_n)=v_n^{(2)}-3v_{n-1}^{(1)}.$$
Let $d\in\{-2,\ldots,p^a\}$. Applying Theorem 2.1(ii) with $h=4$ and $m=5$, we get
$$S_d\eq u_{p^a+d+2}-u_{d+2}\ (\mo\ p)$$
and thus
$$11S_d\eq (v_{p^a+d+2}^{(2)}-v_{d+2}^{(2)})-3(v_{p^a+d+1}^{(1)}-v_{d+1}^{(1)})\ (\mo\ p).$$
Therefore, with the help of Lemma 5.2, we have
$$\align11S_0\eq& (v_{p^a+2}^{(2)}-v_{2}^{(2)})-3(v_{p^a+1}^{(1)}-v_1^{(1)})
\\\eq&\cases2-3(-3)\ (\mo\ p)&\t{if}\ p^a\eq1\ (\mo\ 5),\\-3-3\times2\ (\mo\ p)&\t{if}\ p^a\eq-1\ (\mo\ 5),
\\-4-3(-1)\ (\mo\ p)&\t{if}\ p^a\eq\pm2\ (\mo\ 5);
\endcases
\endalign$$
and
$$\align11S_1
\eq& (v_{p^a+3}^{(2)}-v_{3}^{(2)})-3(v_{p^a+2}^{(1)}-v_2^{(1)})
\\\eq&\cases-18-3(-6)\ (\mo\ p)&\t{if}\ p^a\eq1\ (\mo\ 5),\\16-3\times7\ (\mo\ p)&\t{if}\ p^a\eq-1\ (\mo\ 5),
\\-8-3\times2\ (\mo\ p)&\t{if}\ p^a\eq2\ (\mo\ 5),
\\-5-3\times3\ (\mo\ p)&\t{if}\ p^a\eq-2\ (\mo\ 5).\endcases
\endalign$$
Also,
$$\align11S_{-1}\eq& (v_{p^a+1}^{(2)}-v_{1}^{(2)})-3(v_{p^a}^{(1)}-v_0^{(1)})=v_{p^a+1}^{(2)}-3v_{p^a}^{(1)}
\\\eq&\cases3-3\times1\ (\mo\ p)&\t{if}\ p^a\eq1\ (\mo\ 5),\\-3-3\times0\ (\mo\ p)&\t{if}\ p^a\eq-1\ (\mo\ 5),
\\1-3(-2)\ (\mo\ p)&\t{if}\ p^a\eq2\ (\mo\ 5),
\\-1-3\times1\ (\mo\ p)&\t{if}\ p^a\eq-2\ (\mo\ p).
\endcases
\endalign$$
and
$$\align11S_{-2}\eq& (v_{p^a}^{(2)}-v_{0}^{(2)})-3(v_{p^a-1}^{(1)}-v_{-1}^{(1)})=v_{p^a}^{(2)}-3v_{p^a-1}^{(1)}
\\\eq&\cases0-3\times0\ (\mo\ p)&\t{if}\ p^a\eq1\ (\mo\ 5),\\-1-3\times0\ (\mo\ p)&\t{if}\ p^a\eq-1\ (\mo\ 5),
\\-1-3\times5\ (\mo\ p)&\t{if}\ p^a\eq2\ (\mo\ 5),
\\2-3(-5)\ (\mo\ p)&\t{if}\ p^a\eq-2\ (\mo\ p).
\endcases
\endalign$$
This proves part (i).

(ii) Part (ii) follows from the first congruence in Theorem 2.1(i) with $h=3$ and $m=5$.

(iii) As $C_k^{(3)}=\bi{4k}k-3\bi{4k}{k-1}$ and $\bar C_k^{(3)}=3\bi{4k}k-3\bi{4k}{k+1}$
 for any $k\in\N$, if $p\not=11$ then we can obtain the last two congruences
in Theorem 1.8 by using the congruences on $S_0,S_{\pm1}$ mod $p$ in part (i).

  Below we handle the case $p=11$. This time we turn our resort to Theorem 2.1(i). By (2.4) in the case $h=3$ and $m=5$,
$$\align&\sum_{k=0}^{p^a-1}\f{C_k^{(3)}}{5^k}=\sum_{k=0}^{p^a-1}\f{\bi{4k}k}{5^k}-3\sum_{k=0}^{p^a-1}\f{\bi{4k}{k-1}}{5^k}
\\\eq&-\sum_{r+1}^4\bi 4{r+1}(u_{2-rp^a}-3u_{2-1-rp^a})=-\sum_{r=1}^3\bi 4{r+1}v_{-rp^a}^{(1)}\ (\mo\ p).
\endalign$$
and
$$\align&\sum_{k=0}^{p^a-1}\f{\bar C_k^{(3)}}{5^k}=3\sum_{k=0}^{p^a-1}\f{\bi{4k}k}{5^k}
-\sum_{k=0}^{p^a-1}\f{\bi{4k}{k+1}}{5^k}
\\\eq&-\sum_{r+1}^4\bi 4{r+1}(3u_{2-rp^a}-u_{2+1-rp^a})=\sum_{r=1}^3\bi 4{r+1}v_{1-rp^a}^{(1)}\ (\mo\ p).
\endalign$$
By the proof of Lemma 5.1, $v_n^{(1)}=V_n^{(1)}$ for all $n\in\Z$. Since $p^a=11^a\eq1\ (\mo\ 5)$,
if $\zeta$ is a 5th root of unity
then
$$(1+\zeta)^{-2rp^a}\eq(1+\zeta^{p^a})^{-2r}=(1+\zeta)^{-2r}\ (\mo\ p).$$
Thus
$$v_{-rp^a}^{(1)}=V_{-rp^a}^{(1)}\eq V_{-r}^{(1)}=v_{-rp^a}^{(1)}\ \ (\mo\ p)$$
and
$$v_{1-rp^a}^{(1)}=V_{1-rp^a}^{(1)}\eq V_{1-r}^{(1)}=v_{1-rp^a}^{(1)}\ \ (\mo\ p).$$
Therefore
$$\align\sum_{k=0}^{p^a-1}\f{C_k^{(3)}}{5^k}
\eq&-\sum_{r=1}^3\bi 4{r+1}v_{-r}^{(1)}=-(6v^{(1)}_{-1}+4v_{-2}^{(1)}+v_{-3}^{(1)})
\\\eq&v_1^{(1)}-v_0^{(1)}=u_3-3u_2-(u_2-3u_1)=1\ (\mo\ p).
\endalign$$
and
$$\align\sum_{k=0}^{p^a-1}\f{\bar C_k^{(3)}}{5^k}
\eq&\sum_{r=1}^3\bi 4{r+1}v_{1-r}^{(1)}=6v^{(1)}_{0}+4v_{-1}^{(1)}+v_{-2}^{(1)}
\\\eq&v_1^{(1)}-v_2^{(1)}=u_3-3u_2-(u_4-3u_3)=3\ (\mo\ p).
\endalign$$

 In view of the above, we have completed the proof of Theorem 1.8. \qed

 \heading{6. Proof of Theorem 1.9}\endheading

\medskip
\noindent{\it Proof of Theorem 1.9}. Let $U_0=U_1=U_2=0$, $U_3=1$ and
$$U_{n+4}+\l(4-\f{4^4}{3^3}\r)U_{n+3}+6U_{n+2}+4U_{n+1}+U_n=0\ \quad\t{for}\ n\in\Z.$$
Observe that
$$(1+x)^4-\f{4^4}{3^3}x^3=(x-3)^2\l(x-\f{\al}{27}\r)\l(x-\f{\beta}{27}\r),$$
where $\al+\beta=-14$ and $\al\beta=81$. Let $u_n=u_n(-14,81)$ and $v_n=v_n(-14,81)$ for $n\in\Z$.
By induction,
$$2^5U_n=(6n-11)3^{n-1}+3^{-3(n-1)}(5u_n-11u_{n-1})\ \quad\t{for}\ n\in\Z.$$
This, together with Fermat's little theorem and Theorem 2.1(i) with $h=3$ and $m=4^4/3^3$, yields that
if $d\in\{-2,\ldots,p^a\}$ then
$$\align&-\sum_{k=0}^{p^a-1}\f{3^{3k}}{4^{4k}}\bi{4k}{k+d}
\\\eq&6U_{2+d-p^a}+4U_{2+d-2p^a}+U_{2+d-3p^a}
\\\eq&\f{67}{64}(6d+1)3^{d-2}+\f{5(2u_{d+2-p^a}+36u_{d+2-2p^a}+3^5u_{d+2-3p^a})}{64\times 3^{2d-1}}
\\&-\f{11(2u_{d+1-p^a}+36u_{d+1-2p^a}+3^5u_{d+1-3p^a})}{64\times 3^{2d-1}}
\ (\mo\ p).
\endalign$$

Let $n$ be any integer. Note that $v_n=2u_{n+1}+14u_n$ and $\Delta:=(-14)^2-4\times81=-2^7$. Applying Lemma 3.1 we get
$$u_{n-p^a}\eq-\f7{81}\l(1+\l(\f{-2}{p^a}\r)\r)u_n-\l(\f{-2}{p^a}\r)\f{u_{n+1}}{81}\ (\mo\ p).$$
It follows that
$$u_{n-2p^a}=u_{(n-p^a)-p^a}\eq\f{17+98(\f{-2}{p^a})}{81^2}u_n+\f{14}{81^2}\l(\f{-2}{p^a}\r)u_{n+1}\ (\mo\ p)$$
and
$$u_{n-3p^a}=u_{(n-p^a)-2p^a}\eq\f{(329-805(\f{-2}{p^a}))u_n-115(\f{-2}{p^a})u_{n+1}}{81^3}\ (\mo\ p).$$

Combining the above, for any $d=-2,\ldots,p^a$ we obtain the congruence
$$\aligned&64\sum_{k=0}^{p^a-1}\f{3^{3k}}{4^{4k}}\bi{4k}{k+d}+67(6d+1)3^{d-2}
\\\eq&\f{(1705-482(\f{-2}{p^a}))u_{d+1}-(775+46(\f{-2}{p^a}))u_{d+2}}{27^{d+2}}\ (\mo\ p).
\endaligned\tag6.1$$

Putting $d=0,-1$ in (6.1) we get
$$\sum_{k=0}^{p^a-1}\f{3^{3k}}{4^{4k}}\bi{4k}k\eq\f{44+(\f{-2}{p^a})}{288}\ (\mo\ p)$$
and
$$3\sum_{k=0}^{p^a-1}\f{3^{3k}}{4^{4k}}\bi{4k}{k-1}\eq-\f{220+23(\f{-2}{p^a})}{288}\ (\mo\ p).$$
It follows that
$$\sum_{k=0}^{p^a-1}\f{3^{3k}}{4^{4k}}C_k^{(3)}=\sum_{k=0}^{p^a-1}\f{3^{3k}}{4^{4k}}\(\bi{4k}k-3\bi{4k}{k-1}\)
\eq\f{(\f{-2}{p^a})-1}{12}\ (\mo\ p).$$

By Lemma 3.1,
$$2u_{p^a+1}\eq-14u_1+\l(\f{\Delta}{p^a}\r)v_1=-14-14\l(\f{-2}{p^a}\r)\ (\mo\ p)$$
and
$$2u_{p^a+2}\eq-14u_2+\l(\f{\Delta}{p^a}\r)v_2=196+34\l(\f{-2}{p^a}\r)\ (\mo\ p).$$
Thus, by taking $d=p^a$ in (6.1) we obtain the second congruence in Theorem 1.9.
We are done. \qed

\enddocument